\documentclass[11pt]{amsart}
\usepackage{amssymb}

\newcommand{\Z}{{\mathbb{Z}}}
\newcommand{\Q}{{\mathbb{Q}}}
\newcommand{\F}{{\mathbb{F}}}

\newcommand{\C}{{\mathbb{C}}}
\newcommand{\ba}{{\mathbf{a}}}
\newcommand{\bc}{{\mathbf{c}}}
\newcommand{\bd}{{\mathbf{d}}}
\newcommand{\cO}{{\mathcal{O}}}
\newcommand{\cB}{{\mathcal{B}}}
\newcommand{\cE}{{\mathcal{E}}}
\newcommand{\cH}{{\mathcal{H}}}
\newcommand{\fB}{{\mathfrak{B}}}

\newcommand\Uch{\operatorname{Uch}}
\newcommand\UBr{\operatorname{UBr}}
\newcommand{\Irr}{{\operatorname{Irr}}}
\newcommand{\IBr}{{\operatorname{IBr}}}
\newcommand{\End}{{\operatorname{End}}}
\newcommand{\GL}{{\operatorname{GL}}}
\newcommand{\GU}{{\operatorname{GU}}}

\renewcommand{\leq}{\leqslant}
\renewcommand{\geq}{\geqslant}

\newtheorem{thm}{Theorem}[section]
\newtheorem{lem}[thm]{Lemma}
\newtheorem{cor}[thm]{Corollary}
\newtheorem{prop}[thm]{Proposition}
\newtheorem{conj}[thm]{Conjecture}

\theoremstyle{definition}

\theoremstyle{remark}
\newtheorem{rem}[thm]{Remark}

\begin{document}

\title{Modular principal series representations}
\author{Meinolf Geck}
\address{Department of Mathematical Sciences, King's College,
Aberdeen University, Aberdeen AB24 3UE, Scotland, U.K.}

\email{m.geck@maths.abdn.ac.uk}

\date{February, 2006}

\begin{abstract}
We present a classification of the modular principal series
representations of a finite group of Lie type, in non-describing
characteristic. The proofs rely on the recent progress concerning
the determination of the irreducible representations of
Iwahori--Hecke algebras at roots of unity.
\end{abstract}

\maketitle

\pagestyle{myheadings}

\markboth{Geck}{Modular principal series representations}

\maketitle

%%%%%%%%%%%%%%%%%%%%%%%%%%%%%%%%%%%%%%%%%%%%%%%%%%%%%%%%%%%%%%%%%%%%%%%%%%%%%
\section{Introduction} \label{sec:intro}
Recently, there has been considerable progress in classifying the  
irreducible representations of Iwahori--Hecke algebras at roots of unity.
In this paper, we present an application of these results to modular 
Harish--Chandra series for a finite group of Lie type. 

Let $G$ be a connected reductive algebraic group defined over the 
finite field ${\F}_q$ with $q$ elements. Then $G^F=\{g \in G \mid 
F(g)=g\}$ is a finite group of Lie type, where $F \colon G 
\rightarrow G$ is the corresponding Frobenius map. 

Let $\Irr(G^F)$ be the set of complex irreducible characters of $G^F$. 
It is well-known that $\Irr(G^F)$ is partitioned into Harish--Chandra
series; see, for example, Carter \cite[\S 9.2]{Ca2}. Let us consider
the ``principal series'', denoted $\Irr(G^F,B^F)$ where $B\subseteq G$ 
is an $F$-stable Borel subgroup. By definition, we have $\chi \in 
\Irr(G^F,B^F)$ if $\chi$ appears with non-zero multiplicity in the 
character of the permutation module $\C[G^F/B^F]$. Let $W$ be the 
Weyl group of $G$, with respect to an $F$-stable maximal torus of $B$.
Then $F$ induces an automorphism of $W$ which we denote by the same 
symbol. It is a classical fact (see, for example, Curtis--Reiner 
\cite[\S 68B]{CR2}) that we have a natural bijection
\begin{equation*}
\Irr(W^F) \stackrel{\sim}{\longrightarrow} \Irr(G^F,B^F), \qquad
\rho \mapsto [\rho].\tag{$\spadesuit$}
\end{equation*}
(To be more precise, the bijection depends on the choice of a square root
of $q$; see Lusztig \cite{LuBC}, Geck--Pfeiffer \cite{ourbuch}.) It is the 
purpose of this paper to prove an analogous result for $\ell$-modular Brauer 
characters, where $\ell$ is a prime not dividing $q$. We shall see that, if
$\ell$ is not too small (this will be specified later), then the 
following hold:
\begin{itemize}
\item[1)] the ``$\ell$-modular principal series'' is naturally in 
bijection with a {\em subset} of $\Irr(W^F)$;
\item[2)] that subset is ``generic'', in the sense that it only depends on 
\[ e=\min\{i \geq 2\mid 1+q+q^2+ \cdots +q^{i-1} \equiv 0 \bmod \ell\},\]
but not on the particular values of $q$ and $\ell$;
\item[3)] the bijection in 1) is ``natural'' in the sense that it fits 
into a more general (conjectural) bijection between unipotent characters 
(in the sense of Deligne--Lusztig) and unipotent Brauer characters.
\end{itemize}
If $G^F=\GL_n(q)$, then 1), 2) and 3) reduce to classical results due to 
Dipper and James \cite{DJ0}, \cite{DJ1}, \cite{QuotHom}. Already for 
$G^F=\GU_n(q)$ (the general unitary group), completely new methods are
needed. The results in this paper form one step towards a complete 
description of $\ell$-modular Harish-Chandra series in general; see the 
articles by Hiss, Malle and the author \cite{ghm1}, \cite{ghm2}, 
\cite{lymgh} for a further discussion of this (as yet open) problem. 

To state our results more precisely, we need to introduce some notation.

Every $\chi \in \Irr(G^F)$ has a well-defined {\em unipotent support},
denoted $C_\chi$. This is the unique $F$-stable unipotent conjugacy 
class in $G$, of maximal possible dimension such that a certain average
value of $\chi$ on the $F$-fixed points in that class is non-zero. The 
existence of $C_\chi$ was proved by Lusztig \cite{L6}, assuming that $q$ is 
a sufficiently large power of a sufficiently large prime. These conditions 
on $q$ were removed by Geck--Malle \cite{GeMa}. We set 
\[ \bd_\chi:=\dim \fB_u \qquad (u \in C_\chi)\]
where $\fB_u$ is the variety of Borel subgroups containing $u$.

Now let $\ell$ be a prime not dividing $q$. We shall freely use the general 
notions of the $\ell$-modular representation theory of finite groups (see
Curtis--Reiner \cite[\S 18]{CR2}). Let $\IBr_{\ell}(G^F)$ be the set of 
irreducible Brauer characters of $G^F$. We have equations
\[ \hat{\chi}=\sum_{\varphi \in \IBr_{\ell}(G^F)} d_{\chi,\varphi}\,
\varphi\qquad \mbox{for any  $\chi \in \Irr(G^F)$},\]
where $\hat{\chi}$ denotes the restriction of $\chi$ to the set of 
$\ell$-regular elements of $G^F$ and where  $d_{\chi,\varphi}\in 
{\Z}_{\geq 0}$ are Brauer's $\ell$-modular decomposition numbers. 

Hiss \cite{Hiss1} has shown how to define $\ell$-modular Harish--Chandra
series in general. In this context, the ``$\ell$-modular principal series'', 
denoted $\IBr_{\ell} (G^F,B^F)$, is the set of all Brauer characters which 
are afforded by a simple quotient (or, equivalently, a simple submodule) of
the permutation module $k[G^F/B^F]$, where $k$ is a 
sufficiently large field of characteristic~$\ell$. (See also Dipper 
\cite{QuotHom} for this special case.)  For any $\varphi \in 
\IBr_\ell(G^F,B^F)$, we set 
\[ {\bd}_\varphi^\prime:=\min\{\bd_\chi \mid \chi \in \Irr(G^F,B^F) 
\mbox{ and } d_{\chi,\varphi}\neq 0\}.\]
We can now state the main results of this paper. We will assume 
throughout that $G$ is simple modulo its center.  

\begin{thm} \label{thmA} Assume that $\ell$ is not ``too small'' (see 
Remark~\ref{remA} below). For any $\varphi \in \IBr_{\ell} (G^F,B^F)$, 
there is a unique $\rho=\rho_\varphi \in \Irr(W^F)$ such that 
\[ d_{[\rho],\varphi}\neq 0 \qquad \mbox{and}\qquad 
{\bd}_\varphi^\prime=\bd_{[\rho]}.\]
This yields an injective map $\IBr_\ell(G^F,B^F) \hookrightarrow \Irr(W^F)$, 
$\varphi \mapsto \rho_\varphi$. 
\end{thm}

The proof will be given in \S \ref{sub-fin}; it essentially relies on
the theory of Iwahori--Hecke algebras. These come into play via the
fact, due to Dipper \cite[Cor.~4.10]{QuotHom}, that we have a canonical 
bijection between $\IBr_{\ell}(G^F,B^F)$ and the simple modules of the 
endomorphism algebra
\[ \cE_k:=\End_{kG^F}(k[G^F/B^F]).\]
Using this fact, the desired injection of $\IBr_\ell(G^F,B^F)$ into 
$\Irr(W^F)$ follows from a suitable classification of the simple
$\cE_k$-modules, which is provided by the methods developed in \cite{mykl}, 
\cite{my00}, \cite{GeRo2}; see also the survey in \cite{mylaus}. 

We point out that the statements of Theorem~\ref{thmA} are false when 
$\ell$ is very small, for example, when $\ell=2$ and $G$ is a group of 
type $B_n$, $C_n$ or $D_n$.

\begin{thm} \label{thmB} Keep the hypotheses of Theorem~\ref{thmA}. Then 
the image of the map $\IBr_\ell(G^F,B^F) \hookrightarrow \Irr(W^F)$ only 
depends on $e$ (as defined above).
\end{thm}

Explicit descriptions of the image of the map $\IBr_\ell(G^F,B^F) 
\hookrightarrow \Irr(W^F)$ are now known, through the work of several 
authors: Ariki \cite{Ar1}, \cite{Ar3}, Ariki--Mathas \cite{ArMa}, 
Bremke \cite{Br}, Dipper--James \cite{DJ0}, Dipper--James--Murphy 
\cite{DJM3}, Geck \cite{habil}, \cite{mykl}, Geck--Jacon \cite{GeJa}, 
Geck--Lux \cite{GeLu}, Jacon \cite{Jac0}, \cite{Jac1}, \cite{Jac2}, 
M\"uller \cite{muell}; see Theorem~\ref{thmD}.  Parts of these explicit 
results are needed in the proof of Theorem~\ref{thmB}; see \S 3 for 
the details.

Finally, let us explain why we consider the parametrization in 
Theorem~\ref{thmA} to be ``natural''. Let $\Uch(G^F)\subseteq \Irr(G^F)$ 
be the set of unipotent irreducible characters, as defined by 
Deligne--Lusztig (see \cite{Ca2}, \cite{Lu4}). Let $\UBr_\ell(G^F)$ be 
the set of all $\varphi \in \IBr_\ell(G^F)$ such that $d_{\chi,\varphi} 
\neq 0$ for some $\chi \in \Uch(G^F)$. For $\varphi \in \UBr_\ell(G^F)$, we 
set
\[ {\bd}_\varphi^\prime:=\min\{\bd_\chi\mid \chi \in \Uch(G^F) \mbox{ and }
d_{\chi,\varphi}\neq 0\}.\]

\begin{conj}[Geck \protect{\cite[\S 2.5]{myphd}}, Geck--Hiss 
\protect{\cite[\S 3]{lymgh}}] \label{conj1} Assume that $\ell$ is not 
``too small'' and that the center of $G$ is connected. For any 
$\varphi \in \UBr_{\ell}(G^F)$, there is a unique $\chi=\chi_\varphi 
\in \Uch(G^F)$ such that 
\[ d_{\chi,\varphi}\neq 0 \qquad \mbox{and}\qquad {\bd}_\varphi^\prime=
\bd_{\chi}.\] 
This yields a bijection $\UBr_\ell(G^F) \stackrel{\sim}{\longrightarrow}
\Uch(G^F)$, $\varphi \mapsto \chi_\varphi$. 
\end{conj}

Let $D_{\text{uni}}$ be the matrix of decomposition numbers $d_{\chi,
\varphi}$ where $\chi \in \Uch(G^F)$ and $\varphi \in \UBr_{\ell}(G^F)$.
By \cite{bs1}, we already know that 
\[|\Uch(G^F)|=|\UBr_\ell(G^F)| \qquad \mbox{and}\qquad 
\det D_{\text{uni}}=\pm 1.\]
The above conjecture predicts that $D$ has a 
block triangular shape as follows:
\[ D_{\text{uni}}=\left( \begin{array}{cccc} D_1 & 0 &  \ldots  & 0 \\ 
*  & D_2  & \ddots  & \vdots  \\ \vdots & \ddots & \ddots & 0 \\  
*   & \ldots  &  * & D_r \end{array}\right)\]  
where the rows are partitioned into blocks corresponding to the equivalence
relation ``$\chi \sim \chi' \Leftrightarrow C_\chi=C_{\chi'}$'', the 
blocks are ordered by increasing value of $\bd_\chi$, and the columns 
are ordered via the bijection $\UBr_\ell(G^F) \stackrel{\sim}{\rightarrow} 
\Uch(G^F)$. Furthermore, each $D_i$ is in fact  an identity matrix.

The conjecture is known to be true for $G^F=\GL_n(q)$ (where it follows
from the results of Dipper--James \cite{DJ1}), $G^F=\GU_n(q)$ (see 
\cite{myuni}), and  some explicitly worked-out examples of small rank, 
like $G^F=G_2(q)$ or ${^3\!}D_4(q)$ (see Hiss \cite{H2} and the 
references there).

\begin{cor} \label{cons1} Assume that $\ell$ is not ``too small'' and
that Conjecture~\ref{conj1} holds for $G^F$. Then we have a commutative 
diagram
\begin{center}
\begin{picture}(200,50)
\put( 10, 5){$\UBr_{\ell}(G^F)$}
\put(130, 5){$\Uch(G^F)$}
\put( 75, 7){\vector(1,0){35}}
\put( 87,10){$\sim$}
\put(  0,35){$\IBr_{\ell}(G^F,B^F)$}
\put(120,35){$\Irr(G^F,B^F)$}
\put( 77,37){\vector(1,0){33}}
\put( 77,39){\oval(8,4)[l]}
\put( 25,20){$\cap$}
\put( 33,20){\line(0,1){6}}
\put(145,20){$\cap$}
\put(153,20){\line(0,1){6}}
\end{picture}
\end{center}
where the top horizontal arrow is given by the injection in 
Theorem~\ref{thmA} (composed with the bijection $(\spadesuit)$)
and the bottom horizontal arrow is given by the bijection in 
Conjecture~\ref{conj1}.
\end{cor}

This is an immediate consequence of the results in Section~2 (see the 
argument in the proof of Lemma~\ref{lem1}).

\begin{rem} \label{remA} We will show that the statements in 
Theorems~\ref{thmA} and \ref{thmB} hold under some conditions on $\ell$, 
which are specificied as follows:
\begin{itemize}
\item[(a)] $\ell$ good for $G^F$, and $F=$ identity on $W$; see
Corollary~\ref{eins} and \S\ref{sub-fin}. 
\item[(b)] $\ell$ sufficiently large (explicit bound unknown), any $G$, $F$; 
see \S\ref{sub-fin}.
\item[(c)] $G^F=\GU_n(q)$, any $\ell$; see \S \ref{sub-uni}. There,
we also give an explicit description of the image of the map $\Irr(W^F)
\stackrel{\sim}{\rightarrow}\Irr(G^F,B^F)$.
\end{itemize}
Here, the conditions for $\ell$ to be good are as follows:
\begin{center} 
$\begin{array}{rl} A_n, {^2\!A}_n: & \mbox{no condition}, \\
B_n, C_n, D_n, {^2\!D}_n: & \ell \neq 2, \\
G_2, {^3\!D}_4, F_4, {^2\!E}_6, E_6, E_7: &  \ell \neq 2,3, \\
E_8: & \ell \neq 2,3,5.  \end{array}$
\end{center}
We conjecture that, in general, it is sufficient to assume that $\ell$ is
good for $G^F$. This is supported by the fact that the arguments 
for proving Theorems~\ref{thmA} and \ref{thmB} work whenever $\ell$ is 
good and Lusztig's conjectures on Hecke algebras with unequal parameters 
(as stated in \cite[\S 14.2]{Lusztig03}) hold for $W^F$ and the weight 
function given by the restriction of the length function on $W$ to $W^F$. 
These conjectures are known to hold when $F$ is the identity on $W$ 
(the ``equal parameter case''); see \cite[Chap.~15]{Lusztig03}. A sketch 
proof for the case where $F$ is not the identity is given in 
\cite[Chap.~16]{Lusztig03}.

We also expect that Conjecture~\ref{conj1} holds whenever $\ell$ is good 
for $G^F$.
\end{rem}

\begin{rem} \label{remB} As already remarked above, there is a natural 
bijection between $\IBr_\ell(G^F,B^F)$ and the set of simple modules of 
the algebra $\cE_k$. Assume now that $W^F$ is of type $B_n$ or $D_n$. 
As we shall see, in these cases, the image of the map in 
Theorem~\ref{thmA} is explicitly given by the results of Jacon 
\cite{Jac0}, \cite{Jac1}, \cite{Jac2}, and Jacon and the author 
\cite{GeJa}. 

Note that Ariki \cite{Ar2} (see also Ariki--Mathas \cite{ArMa}) gives 
a different parametri\-zation of the simple $\cE_k$-modules by a subset 
of $\Irr(W^F)$, based on the Dipper--James--Murphy theory \cite{DJM3} 
of Specht modules.
\end{rem}

%%%%%%%%%%%%%%%%%%%%%%%%%%%%%%%%%%%%%%%%%%%%%%%%%%%%%%%%%%%%%%%%%%%%%%%%%%%%%
\section{Decomposition numbers} \label{sec2}
We keep the basic set-up of Section~1. The purpose of this section is to
show how the proof of Theorem~\ref{thmA} can be reduced to an analogous
statement about decomposition numbers of Iwahori--Hecke algebras.  This 
involves three major ingredients:
\begin{itemize}
\item Iwahori's realization of the endomorphism algebra $\cE_k$ as a 
specialization of a ``generic'' algebra; see Proposition~\ref{thmI}.
\item Dipper's interpretation of the decomposition numbers 
$d_{\chi,\varphi}$ for $\chi \in \Uch(G^F)$ and $\varphi\in 
\UBr_\ell(G^F,B^F)$ in terms of $\cE_k$; see Proposition~\ref{dip}.
\item Lusztig's Hecke algebra interpretation for the invariants 
$\bd_{[\rho]}$ where $\rho \in \Irr(W^F)$; see Theorem~\ref{uni2}.
\end{itemize}
Once this reduction is achieved, we can apply the results of Jacon, 
Rouquier and the author on the existence of so-called ``canonical
basic sets'' for Iwahori--Hecke algebras. (For a survey, see 
\cite{mylaus}).

Let $K$ be a sufficiently large finite field extension of $\Q$; specifically, 
we require that $K$ contains all $|G^F|$th roots of unity and an element
$t$ such that $t^2=q$. As is well-known, all complex irreducible characters 
of $G^F$ can be realized over $K$, so we can actually regard $\Irr(G^F)$ 
as the set of irreducible $K$-characters of $G^F$.  A similar remark applies
to $\Irr(W^F)$.

Let $\cO$ be a discrete valuation ring in $K$, with residue field $k$ of 
characteristic $\ell>0$. Then it is also known that $k$ is a splitting for
$G^F$. Having chosen $\cO$, the $\ell$-modular Brauer character of a 
simple $kG^F$-module is well-defined. 

For $R \in \{K,\cO,k\}$, we define  
\[ \cE_R:=\End_{RG^F}(R[G^F/B^F])\]
where $R[G^F/B^F]$ is the $RG^F$-permutation module on the cosets of $B^F$. 
By Iwahori's theorem (see \cite[\S 8.4]{ourbuch} or \cite[\S 10.10]{Ca2}), 
the algebra $\cE_R$ has a standard basis indexed by the elements of $W^F$; 
furthermore, $\cE_R$ can be obtained from a ``generic algebra''. This is 
done as follows.

\subsection{The generic Iwahori--Hecke algebra} \label{sub-gen}

The group $W$ is a finite Coxeter group with generating set $S$. Let 
$l \colon W \rightarrow \Z_{\geq 0}$ be the corresponding length function. 
We have $F(S)=S$. Let $\bar{S}$ be the set of $F$-orbits on $S$. Then 
$W^F$ is a finite  Coxeter group with generating set $\{w_J \mid J \in 
\bar{S}\}$, where $w_J$ is the longest element in the parabolic subgroup 
generated by $J$.  Let $\bar{l}\colon W^F \rightarrow \Z_{\geq 0}$ be the 
corresponding length function. Furthermore, let $L\colon W^F \rightarrow 
{\Z}_{\geq 0}$ be the restriction of the length function on $W$ to $W^F$. 
Then $L$ is a weight function in the sense of Lusztig; see 
\cite[Lemma~16.2]{Lusztig03}.

Let $A={\cO}[v,v^{-1}]$ be the ring of Laurent polynomials over $\cO$ in an
indeterminate~$v$. Let $\cH$ be the generic Iwahori--Hecke algebra of $W^F$ 
over $A$. Thus, $\cH$ is free as an $A$-module, with a basis $\{T_w \mid w 
\in W^F\}$, such that the following relations hold:
\begin{align*}
T_{ww'} &= T_w T_{w'} \quad \mbox{whenever $\bar{l}(ww')=\bar{l}(w)+
\bar{l} (w')$},\\ T_{w_J}^2 & =u^{L(w_J)}T_1+(u^{L(w_J)}-1)T_{w_J} 
\quad  \mbox{for any $J\in \bar{S}$, where  $u=v^2$}.
\end{align*}
Now there is a canonical ring homomorphism $\theta_0 \colon A \rightarrow \cO$
such that $\theta_0(v)=t$. Similarly, there is a canonical ring homomorphism 
$\theta \colon A \rightarrow k$ such that $\theta(v)=\bar{t}$, where 
$\bar{t}$ denotes the image of $t$ in $k$. Thus, $\theta$ is the
composition of $\theta_0$ with the canonical map from $\cO$ onto~$k$. 
For $R \in \{K,\cO,k\}$, we write $\cH_R=R \otimes_A \cH$ where $R$ is
regarded as an $A$-module via the ring homomorphism $\theta_0$ or $\theta$.

\begin{prop}[Iwahori] \label{thmI} We have $\cH_R\cong \cE_R$.
\end{prop}

Let $K(v)$ be the field of fractions of $A$ and write $\cH_{K(v)}=
K(v) \otimes_A \cH$. It is known that the algebra $\cH_{K(v)}$ is split 
semisimple (see Lusztig \cite{LuBC} for the case where $L$ is
constant on $S$ and \cite[Chap.~9]{ourbuch} in general). By Tits'
Deformation Theorem (see \cite[\S 8.1]{ourbuch}), the ring homomorphism 
$A \rightarrow K$, $v \mapsto 1$, induces a bijection between $\Irr(W^F)$
and the set of simple $\cH_{K(v)}$-modules, up to isomorphism. Let us
write 
\[ \Irr(W^F)=\{\rho^\lambda \mid \lambda \in \Lambda\}\]
where $\Lambda$ is a finite indexing set. Then, for each $\lambda \in
\Lambda$, there is a well-defined simple $\cH_{K(v)}$-module $E_v^\lambda$.
We have $\mbox{trace}(T_w,E_v^\lambda)\in A$ for all $w\in W^F$, and 
$E^\lambda_v$ is uniquely determined by the condition that 
\[ \rho^\lambda(w)=\mbox{trace}(T_w,E_v^\lambda) |_{v=1} \qquad 
\mbox{for all $w\in W^F$}.\] 
Similarly, for each $\lambda \in \Lambda$, there is a corresponding
simple $\cH_K$-module $E^\lambda_{t}$, which is uniquely determined
by the condition that  
\[ \mbox{trace}(T_w,E^\lambda_{t})=\mbox{trace}(T_w,E_v^\lambda) |_{v=t}
\qquad \mbox{for all $w\in W^F$}.\] 
In this context, the ``classical'' bijection (mentioned in Section~1)
\begin{equation*}
\Irr(W^F) \stackrel{\sim}{\longrightarrow} \Irr(G^F,B^F), \qquad 
\rho \mapsto [\rho],\tag{$\spadesuit$}
\end{equation*}
is obtained as follows. Let $\rho\in \Irr(W^F)$ and write $\rho=
\rho^\lambda$ where $\lambda \in \Lambda$. Now regard $E^\lambda_{t}$ 
as a simple $\cE_K$-module via the isomorphism in Proposition~\ref{thmI}. 
Then, by standard results on endomorphism algebras (see, for example, 
\cite[Prop.~8.4.4]{ourbuch} or \cite[\S 68B]{CR2}), the module 
$E^\lambda_{t}$ corresponds to a well-defined irreducible character 
in $\Irr(G^F,B^F)$, which is denoted $[\rho]$.

\subsection{The decomposition matrix of $\cH$} \label{sub-hec}

The ring homomorphism $\theta\colon A \rightarrow k$ induces a 
decomposition map
\[ d_\theta \colon R_0(\cH_{K(v)})\rightarrow R_0(\cH_k)\]
from the Grothendieck group of finitely generated $\cH_{K(v)}$-modules 
to the Grothendieck group of finitely generated $\cH_k$-modules (see 
\cite[\S 2]{mybourb} or \cite[\S 4]{mylaus}). Given  a simple 
$\cH_{K(v)}$-module $E$ and a simple $\cH_k$-module $M$, we denote 
\[ (E:M)= \begin{array}{l} \mbox{multiplicity of $M$ in the image of 
$E$}\\\mbox{(under the decomposition map $d_\theta$)}.
\end{array}\]
Finally, let  us write the set of simple $\cH_k$-modules (up to 
isomorphism) as 
\[ \{ M^\mu \mid \mu \in \Lambda^\circ\} \quad \mbox{where $\Lambda^\circ$
is a finite indexing set}.\]
Using Dipper's Hom functors \cite{QuotHom}, the set $\Lambda^\circ$ also 
canonically parametrizes the $\ell$-modular principal series of $G^F$:
\[ \IBr_\ell(G^F,B^F)=\{ \varphi^\mu \mid \mu \in \Lambda^\circ\}.\]
With these notations, we can now state the following basic result. 

\begin{prop}[Dipper] \label{dip} Let $\chi\in \Uch(G^F)$ and $\mu 
\in \Lambda^\circ$.  Then 
\[ d_{\chi,\varphi^\mu}=\left\{\begin{array}{cl} (E^\lambda_v:M^\mu) &
\quad \mbox{if $\chi=[\rho^\lambda]$ where $\lambda\in \Lambda$},\\
0 & \quad \mbox{otherwise}.\end{array}\right.\]
\end{prop}

In fact, Dipper's original result in \cite[Corollary~4.10]{QuotHom} 
works with the decomposition numbers of the algebra $\cE_{\cO}$. The above
formulation takes into account the factorization result in 
\cite[3.3]{mybourb}, in order to lift the statement to the generic
algebra $\cH$.

\subsection{Lusztig's $\ba$-function} \label{sub-a}

The algebra $\cH$ is symmetric, with symmetrizing trace $\tau \colon
\cH \rightarrow A$ given by $\tau(T_1)=1$ and $\tau(T_w)=0$ for $1 \neq
w \in W^F$. Since $\cH_{K(v)}$ is split semisimple and $A$ is integrally
closed in $K(v)$, there are well-defined Laurent polynomials $\bc_\lambda 
\in A$ such that
\[ \tau(T_w)=\sum_{\lambda \in \Lambda} \frac{1}{\bc_\lambda} \, 
\mbox{trace} (T_w,E^\lambda_v) \qquad \mbox{for all $w\in W^F$}.\]
This follows from a general argument concerning symmetric algebras;
see \cite[Chapter~7]{ourbuch}. Since the weight function $L$ arises from
a finite group of Lie type, we actually have $\bc_\lambda \in {\Z}[u,
u^{-1}]$ and
\[ [\rho^\lambda](1)=\frac{[G^F:B^F]}{\bc_\lambda(q)} \qquad \mbox{for all
$\lambda \in \Lambda$};\]
see Curtis--Reiner \cite[\S 68C]{CR2}. Let us write 
\[ \bc_\lambda=f_\lambda u^{-\ba_\lambda}+\mbox{combination of higher 
powers of $u$},\]
where $f_\lambda\neq 0$ and $\ba_\lambda\geq 0$ are integers. This yields
Lusztig's $\ba$-function
\[ \Irr(W^F) \rightarrow \Z_{\geq 0}, \qquad \rho^\lambda \mapsto 
\ba_\lambda.\]
This function is related to the $\bd$-invariants considered in Section~1 
by the following result.

\begin{thm}[Lusztig, Geck--Malle] \label{uni2}
Let $\lambda \in \Lambda$ and $C$ be the unipotent support  of 
$[\rho^\lambda]\in \Irr(G^F,B^F)$. Then we have
\[ \ba_\lambda=\bd_{[\rho^\lambda]}=\dim \fB_u \qquad (u \in C).\]
\end{thm}

This was proved by Lusztig \cite[\S 10]{L6}, assuming that $q$ is a 
power of a sufficiently large prime. These conditions were removed
by \cite[Prop.~3.6]{GeMa}.

\subsection{Canonical basic sets} \label{sub-can}

Using the $\ba$-function, we now define the notion of a ``canonical 
basic set'', following \cite[Def.~4.13]{mylaus}. For any $\mu \in 
\Lambda^\circ$, we set
\[ {\ba}_\mu^\prime=\min \{\ba_\lambda \mid \lambda \in \Lambda 
\mbox{ and } (E_v^\lambda:M^\mu)\neq 0\}.\]
Let $\iota \colon \Lambda^\circ \rightarrow \Lambda$ be an injective map 
and write $\cB_{k,q}=\iota(\Lambda^\circ)\subseteq \Lambda$. We say that 
$\cB_{k,q}$ is a {\em canonical basic set} for $\cH_k$ if the following 
conditions are satisfied. 
\begin{align*}
(E_v^{\,\iota(\mu)}:M^\mu) & =1\qquad\mbox{for all $\mu\in\Lambda^\circ$},\\
(E_v^{\lambda}:M^\mu) & =0 \qquad\mbox{unless $\ba_\mu^\prime<
\ba_\lambda$ or $\lambda=\iota(\mu)$}.\
\end{align*}
Note that, if it exists, then $\iota$ is uniquely determined by these
conditions: given $\mu \in \Lambda^\circ$, the image $\iota(\mu)$ is the
unique $\lambda \in \Lambda$ such that
\[ (E_v^\lambda:M^\mu)\neq 0 \qquad \mbox{and}\qquad {\ba}_\mu^\prime=
{\ba}_\lambda.\]
The following result provides the desired reduction of Theorem~\ref{thmA} to
a statement entirely in the framework of Iwahori--Hecke algebras.

\begin{lem} \label{lem1} Assume that $\cH_k$ admits a canonical basic set
$\cB_{k,q}=\iota(\Lambda^\circ)$, as above. Then the statements in 
Theorem~\ref{thmA} hold for $G^F$, where
\[ \rho_{\varphi^\mu}=\rho^{\iota(\mu)} \qquad \mbox{for all $\mu \in
\Lambda^\circ$}.\]
Furthermore, assuming that Conjecture~\ref{conj1} holds, we have a 
commutative diagram as in Corollary~\ref{cons1}.
\end{lem}

\begin{proof} This is an immediate consequence of Proposition~\ref{dip}, 
taking into account the identities in Theorem~\ref{uni2}.
\end{proof}

Now we have the following general existence result.

\begin{thm}[Geck \cite{mykl}, \cite{my00}, \cite{mylaus} and
Geck--Rouquier \cite{GeRo2}] \label{gero} Assume that Lusztig's
conjectures {\bf (P1)--(P14)} in \cite[\S 14.2]{Lusztig03} and a certain
weak version of {\bf (P15)} (as specified in \cite[5.2]{mylaus}) hold
for $\cH$. Assume further that $\ell$ is good for $G^F$. Then $\cH_k$ 
admits a canonical basic set.
\end{thm}

\begin{rem} \label{remC}
(a) The above result is proved by a general argument, using
deep properties of the Kazhdan--Lusztig basis of $\cH$. These properties 
are known to hold, for example, when $F$ is the identity on $W$ (the ``equal 
parameter case''); see Lusztig \cite[Chap.~15]{Lusztig03}. A sketch
proof for the case where $F$ is not the identity is given by Lusztig
\cite[Chap.~16]{Lusztig03}.

(b) If $\ell$ is not good, it is easy to produce examples in which
a canonical basic set does not exist; see \cite[4.15]{mylaus}.
\end{rem}

\begin{cor} \label{eins} Assume that $\ell$ is a good prime for $G^F$ and
$F$ is the identity on $W$. Then the statements in Theorem~\ref{thmA}
and Corollary~\ref{cons1} hold for $G^F$.
\end{cor}

\begin{proof} By Lusztig \cite[Chap.~15]{Lusztig03}, the hypotheses of 
Theorem~\ref{gero} are satisfied and so $\cH_k$ admits a  canonical basic 
set. It remains to use Lemma~\ref{lem1}.
\end{proof}

\subsection{The finite unitary groups} \label{sub-uni}

Let $G=\GL_n(\overline{\F}_q)$ and $F$ be such that $G^F=\GU_n(q)$, the
finite general unitary group. Then $W^F$ is a Coxeter group of type $B_m$, 
where $m=n/2$ (if $n$ is even) or $m=(n-1)/2$ (if $n$ is odd). Writing 
$n=2m+s$, the weight function $L \colon W^F \rightarrow \Z_{\geq 0}$ 
is given by 
\begin{center}
\begin{picture}(250,30)
\put(  3, 08){$B_m$}
\put( 40, 08){\circle{10}}
\put( 44, 05){\line(1,0){33}}
\put( 44, 11){\line(1,0){33}}
\put( 81, 08){\circle{10}}
\put( 86, 08){\line(1,0){29}}
\put(120, 08){\circle{10}}
\put(125, 08){\line(1,0){20}}
\put(155, 05){$\cdot$}
\put(165, 05){$\cdot$}
\put(175, 05){$\cdot$}
\put(185, 08){\line(1,0){20}}
\put(210, 08){\circle{10}}
\put( 30, 18){$2s{+}1$}
\put( 77, 18){$2$}
\put(117, 18){$2$}
\put(207, 18){$2$}
\end{picture}
\end{center}
All unipotent classes in $G$ are $F$-stable, and they are naturally 
labelled by the partitions of $n$. Given $\lambda \vdash n$, let 
$C_\lambda$ be the class containing matrices of Jordan type $\lambda$. 
For $u \in C_\lambda$, we have 
\[ \dim \fB_u=n(\lambda):=\sum_{i=1}^r (i-1)\lambda_i\]
where $\lambda=(\lambda_1\geq \lambda_2\geq \ldots \geq \lambda_r>0)$.
The unipotent characters are also labelled by the partitions of $n$. Thus, 
we can write
\[ \Uch(G^F)=\{\chi^\lambda \mid \lambda \vdash n\}.\]
For example, $\chi^{(n)}$ is the unit and $\chi^{(1^n)}$ is
the Steinberg character. We have 
\[ C_{\chi^\lambda}=C_{\lambda} \qquad \mbox{for all $\lambda\vdash n$}.\]
See Lusztig \cite[\S 9]{L1}, Kawanaka \cite{Kaw1} and Carter 
\cite[Chap.~13]{Ca2} for further details. Now let $\ell$ be any prime 
not dividing $q$. By \cite{myuni}, we have a parametrization 
$\UBr_\ell(G^F) =\{\varphi^\mu \mid \mu \vdash n\}$ 
such that the following hold:
\begin{align*}
d_{\chi^\lambda,\varphi^\lambda} &=1 \qquad 
\mbox{for all $\lambda \vdash n$},\\
d_{\chi^\lambda,\varphi^\mu} &=0 \qquad \mbox{unless $\lambda 
\trianglelefteq \mu$},
\end{align*}
where $\trianglelefteq$ denotes the dominance order on partitions. (The 
proof essentially relies on Kawanaka's theory \cite{Kaw1} of 
generalized Gelfand--Graev representations.)

Now, it is known that, for any $\nu,\nu'\vdash n$, we have 
$\nu\trianglelefteq \nu'\Rightarrow n(\nu')\leq n(\nu)$, with equality 
only if $\nu=\nu'$ (see, for example, \cite[Exc.~5.6]{ourbuch}). 
Hence the above conditions on the decomposition numbers can also
be phrased as:
\begin{align*}
d_{\chi^\lambda,\varphi^\lambda} &=1 \qquad 
\mbox{for all $\lambda \vdash n$},\\
d_{\chi^\lambda,\varphi^\mu} &=0 \qquad \mbox{unless $n(\mu)<n(\lambda)$
or $\lambda=\mu$}.
\end{align*}
Using the formula for $\dim \fB_u$, it is now clear that 
Conjecture~\ref{conj1} holds. 

Now let us consider the principal series characters. The set $\Irr(W^F)$ is 
naturally parametrized by the set $\Lambda$ of pairs of partitions of 
total size $m$. The inclusion $\Irr(G^F,B^F)\subseteq \Uch(G^F)$ corresponds 
to an embedding of $\Lambda$ into the set of partitions of $n$, which is
explicitly described in the appendix of \cite{fs1}, based on Lusztig 
\cite[\S 9]{L1}. (The description involves the notions of the $2$-core 
and the $2$-quotient of a partition.)

Using Proposition~\ref{dip} and the identity in Theorem~\ref{uni2}, we 
conclude that the statements in Theorem~\ref{thmA} and Corollary~\ref{cons1} 
also hold for $G^F$. Note, however, that these are pure existence results! 
An explicit combinatorial description of the image of the map $\IBr_\ell
(G^F,B^F)\hookrightarrow \Irr(W^F)$ (or, equivalently, of a canonical basic
set for $\cH_k$) is much harder to obtain. The complete answer was only 
achieved quite recently; see \cite{GeJa}.

If $e=2$, then that image  is given by a class of bipartitions described
in \cite[Theorem~3.4]{GeJa}; if $e>2$ is twice an odd number, then that 
image  is given by a class of bipartitions defined by Foda et al. 
\cite{FLOTW}; see \cite[Theorem~5.4]{GeJa}; otherwise, that image is 
given by the set of all pairs of $e$-regular partitions of total size 
$m$; see  \cite[Theorem~3.1]{GeJa}. (These results hold for any prime 
$\ell$ not dividing $q$.)

%%%%%%%%%%%%%%%%%%%%%%%%%%%%%%%%%%%%%%%%%%%%%%%%%%%%%%%%%%%%%%%%%%%%%%%%%%%%%
\section{Independence of canonical basic sets} \label{sec:indie}

In order to establish the independence statement in Theorem~\ref{thmB}, we 
use a technique originally developed in \cite[\S 4]{mybaum} (see also 
\cite[\S 2]{GeRo1} and \cite[\S 2]{mybourb}), namely, a 
factorization of the decomposition map 
\[ d_\theta \colon R_0(\cH_{K(v)}) \rightarrow R_0(\cH_k).\]
Let $e=\min\{j \geq 2\mid 1+q+q^2+\cdots +q^{j-1} \equiv 0 \bmod \ell\}$ as
before and $\zeta_e:=\exp(2\pi i/e) \in \C$. Choosing $\cO$ suitably, we may 
assume that $\cO$ contains an element $\zeta_{e}'$ such that $\zeta_e'^2=
\zeta_e$ and $t,\zeta_e'$ have the same image in $k$. Then we have a 
canonical ring homomorphism $\theta_e\colon A \rightarrow \cO$ such that 
$\theta_e(v)=\zeta_{e}'$, and $\theta$ is the composition of $\theta_e$ and 
the canonical map $\cO \rightarrow k$. We consider the Iwahori--Hecke algebra 
\[ \cH_{\cO}^{(e)}=\cO \otimes_{A} \cH\]
where $\cO$ is considered as an $A$-module via $\theta_e$. Thus, in 
$\cH_{\cO}^{(e)}$, we have 
\[ T_{w_J}^2  =\zeta_e^{L(w_J)}T_1+(\zeta_e^{L(w_J)}-1)T_{w_J}
\quad  \mbox{for any $J\in \bar{S}$}.\]
The map $\theta_e \colon A\rightarrow \cO$ induces a well-defined 
decomposition map 
\[ d_e \colon R_0({\cH}_{K(v)}) \rightarrow R_0(\cH_K^{(e)})\]
between the Grothendieck groups of $\cH_{K(v)}$ and 
$\cH_K^{(e)}=K\otimes_{\cO} \cH_{\cO}^{(e)}$. Similarly, the canonical
map $\cO \rightarrow k$ induces a decomposition map 
\[ d' \colon R_0(\cH_K^{(e)}) \rightarrow R_0(\cH_k).\]
By \cite[2.6]{mybourb}, we then have the following factorization of 
$d_\theta$:
\[\mbox{\begin{picture}(200,55)
\put(5,45){$R_0({\cH}_{K(v)})$}
\put(65,47){\vector(1,0){80}}
\put(100,52){$\scriptstyle{d_\theta}$}
\put(155,45){$R_0({\cH}_k)$}
\put(48,35){\vector(3,-2){34}}
\put(55,17){$\scriptstyle{d_e}$}
\put(134,13){\vector(3,2){34}}
\put(158,17){$\scriptstyle{d'}$}
\put(88,05){$R_0({\cH}_K^{(e)})$}
\end{picture}}\]
Using this factorisation, we obtain the following result which provides a 
plan for proving Theorem~\ref{thmB}. 

\begin{lem}[Jacon \protect{\cite[The\'or\`eme~3.1.3]{Jac0}}] \label{lem3} 
Assume that the following hold:
\begin{itemize}
\item[(a)] $\cH_k$ admits a canonical basic set, say $\cB_{k,q}$.
\item[(b)] $\cH_K^{(e)}$ admits a canonical basic set, say $\cB_{e}$.
\item[(c)] $\cH_k$ and $\cH_K^{(e)}$ have the same number of simple modules 
(up to isomorphism).
\end{itemize}
Then we have $\cB_e=\cB_{k,q}$. In particular, $\cB_{k,q}$ 
only depends on $e$.
\end{lem}

\begin{proof} In order to illustrate the use of the above factorization,
we give the argument here. Recall our notation $\{E_v^\lambda \mid \lambda 
\in \Lambda\}$ and $\{M^\mu \mid \mu \in \Lambda^\circ\}$ for the simple 
modules of $\cH_{K(v)}$ and $\cH_k$, respectively. 

By (c), we have a labelling of the simple $\cH_K^{(e)}$-modules
\[ \{M_e^\mu \mid \mu \in \Lambda_e^\circ\}\qquad \mbox{where}\qquad
|\Lambda_e^\circ|=|\Lambda^\circ|.\]
The above factorization means that we have the following identity:
\[ (E^\lambda_v :M^\mu)=\sum_{\nu \in \Lambda_e^\circ} (E_v^\lambda:M_e^\nu)
\cdot (M_e^\nu:M^\mu)\]
for all $\lambda \in \Lambda$ and $\mu \in \Lambda^\circ$, where
$(E_v^\lambda:M_e^\nu)$ are the decomposition numbers of $\cH_K^{(e)}$.
By (a), we have a canonical basic set $\cB_{k,q}=\iota(\Lambda^\circ)$ for
$\cH_k$, where $\iota \colon \Lambda^\circ \rightarrow \Lambda$ is an
injection. By (b), we have a canonical basic set $\cB_{e}=\iota_e
(\Lambda_e^\circ)$ for $\cH_K^{(e)}$, where $\iota_e \colon 
\Lambda_e^\circ \rightarrow \Lambda$ is an injection. 

We now define a map $\beta\colon \Lambda^\circ \rightarrow 
\Lambda_e^\circ$, as follows. Let $\mu \in \Lambda^\circ$ and $\lambda=
\iota(\mu)$. Then $(E_v^\lambda:M^\mu)=1$ and so there is a unique 
$\mu_0 \in \Lambda_e^\circ$ such that $(E_v^\lambda:M_e^{\mu_0})\neq 0$ 
and $(M_e^{\mu_0}:M^\mu)\neq 0$; we set $\beta(\mu)=\mu_0$. Now we claim
that 
\begin{equation*}
 \iota(\mu)=\iota_e(\beta(\mu)) \qquad \mbox{for all $\mu \in 
\Lambda^\circ$}.\tag{$*$}
\end{equation*}
This is seen as follows. Let $\mu\in \Lambda^\circ$ and $\lambda=
\iota(\mu)$. By the construction of $\mu_0=\beta(\mu)$, we have 
$(E_v^\lambda:M_e^{\mu_0})=1$ and $(M_e^{\mu_0}:M^\mu)=1$. It remains 
to show that $\ba_{\lambda'}> \ba_{\lambda}$ for any $\lambda' \in 
\Lambda$ such that $\lambda\neq \lambda'$ and $(E_v^{\lambda'}:
M_e^{\mu_0}) \neq 0$. Indeed, we have
\begin{align*}
(E^{\lambda'}_v :M^\mu)&=\sum_{\nu \in \Lambda_e^\circ} 
(E_v^{\lambda'}: M_e^\nu) \cdot (M_e^\nu:M^\mu)\\
&=(E_v^{\lambda'}: M_e^{\mu_0}) \cdot (M_e^{\mu_0}:M^\mu)+
\mbox{further terms},
\end{align*}
and so this decomposition number is non-zero. Since $\cB_{k,q}$ is a 
canonical basic set, we must have $\ba_{\lambda'}>\ba_\lambda$, as desired.
Thus, ($*$) is proved.  This relation implies that $\cB_{k,q} \subseteq 
\cB_e$.  Since $\cB_{k,q}$ and $\cB_e$ have the same cardinality, we 
conclude that $\cB_{k,q}=\cB_e$, as desired.
\end{proof}

We are now going to show that the conditions in Lemma~\ref{lem3} are
satisfied.

\begin{thm} \label{thmD} The algebra $\cH_K^{(e)}$ admits a canonical 
basic set. This canonical basic set is explicitly known in all cases.
\end{thm}

\begin{proof} This is proved by a combination of various results.  
If $F$ acts as the identity on $W$, then the existence of a canonical basic
is covered by Theorem~\ref{gero} (which also works for $K$ instead of $k$). 
For type $A_{n-1}$, an explicit description is given in \cite[Exp.~3.5]{mykl} 
(based on Dipper--James \cite{DJ0}); for type $B_n$ (with equal parameters) 
and type $D_n$, see Jacon \cite{Jac0}, \cite{Jac1}, \cite{Jac2}; for the 
exceptional types $G_2$, $F_4$, $E_6$, $E_7$, $E_8$ (with equal parameters), 
see the tables of Jacon  \cite[\S 3.3]{Jac0}. These tables are derived from
the results on decomposition numbers by Geck--Lux \cite{GeLu}, Geck
\cite{mye6}, \cite{habil}, and M\"uller \cite{muell}.

Assume now that $F$ is not the identity on $W$. Then $W^F$ is of type $G_2$ 
(with parameters $q^3,q$), $F_4$ (with parameters $q^2,q^2,q,q$) or $B_m$.  
In the first two cases, the existence of a canonical basic set follows from 
the explicit determination of the decomposition matrices by 
\cite[Satz~3.13.1]{myphd} (type $G_2$) and Bremke \cite{Br} (type $F_4$). 
The $\ba$-function for type $F_4$ (with possibly unequal parameters) is 
printed in \cite[Exp.~4.8]{mylaus}. 

For the convenience of the reader, let us give the details for type $G_2$ 
with parameters $q^3,q$. Denote the two Coxeter generators of $W^F$ by 
$\alpha$ and $\beta$, where $L(\alpha)=3$ and $L(\beta)=1$. By 
\cite[8.3.1]{ourbuch}, the algebra $\cH_{K(v)}$ has four $1$-dimensional 
irreducible representations:
\[ \begin{array}{rlll}
\text{ind} \colon & T_{\alpha} \mapsto u^3, & &T_{\beta} \mapsto u; \\
\varepsilon \colon & T_{\alpha} \mapsto -1, & &T_{\beta} \mapsto -1; \\
\varepsilon_1 \colon & T_{\alpha} \mapsto u^3, && T_{\beta} \mapsto -1; \\
\varepsilon_2 \colon & T_{\alpha} \mapsto -1,& & T_{\beta} \mapsto u;
\end{array}\]
and two $2$-dimensional representations:
\[ \rho_\delta \colon  T_{\alpha} \mapsto \left( \begin{array}{cc} 
-1 & 0 \\ u^2+\delta u+1 & u^3 \end{array} \right), \qquad T_\beta 
\mapsto \left( \begin{array}{cc} u & u\\ 0 & -1 \end{array} \right),\]
where $\delta=\pm 1$. The $\ba$-invariants are obtained from
\cite[8.3.4]{ourbuch}. The decomposition numbers are given as follows:
\[ \begin{array}{|ccc|cccc|ccccc|cccc|ccccc|} \hline \rho_\lambda & 
\ba_\lambda& & \multicolumn{3}{c}{e=2} & & \multicolumn{4}{c}{e=3} && 
\multicolumn{3}{c}{e=6} & & \multicolumn{5}{c|}{e=12} \\ \hline
\text{ind}&   0&&1&.&.&&1&.&.&.&& 1 & . & . && 1 & . & . & . & . \\
\varepsilon_1&1&&1&.&.&&.&1&.&.&& . & 1 & . && . & 1 & . & . & . \\
\rho_+ &      3&&.&1&.&&.&1&1&.&& . & . & 1 && 1 & . & 1 & . & . \\
\rho_- &      3&&.&.&1&&1&.&.&1&& 1 & 1 & . && . & . & . & 1 & . \\
\varepsilon_2&7&&1&.&.&&.&.&1&.&& 1 & . & . && . & . & . & . & 1 \\
\varepsilon& 12&&1&.&.&&.&.&.&1&& . & 1 & . && . & . & 1 & . & . \\
\hline \end{array}\] 
This yields the following canonical basic sets:
\begin{align*}
e=2 &: \quad \{\mbox{ind},\;\rho_+,\;\rho_-\},\\
e=3 &: \quad \{\mbox{ind},\;\varepsilon_1,\;\rho_+,\;\rho_-\},\\
e=6 &: \quad \{\mbox{ind},\;\varepsilon_1,\;\rho_+\},\\
e=12 &: \quad \{\mbox{ind},\;\varepsilon_1,\;\rho_+,\;\rho_-,\;
\varepsilon_2\}.
\end{align*}
For all other values of $e$, the decomposition matrix is the identity
matrix. Thus, in each case, we see that a canonical basic sets exists. 

Finally, the existence of a canonical basic set in type $B_m$ (for any 
choice of the parameters) is established by Geck--Jacon \cite{GeJa}.
\end{proof}

In those cases where the hypotheses of Theorem~\ref{gero} (concerning
Lusztig's conjectures on Hecke algebras with unequal parameters) are
not yet known to hold, we now have the following result:

\begin{cor} \label{lem4} There is a constant $N$, depending only on $W$, 
such that $\cH_k$ admits a canonical basic set, whenever $\ell>N$.
\end{cor}

\begin{proof} By \cite[Prop.~5.5]{mybaum} (see also \cite[\S 2.7]{mybourb}),
there is a constant $N$, depending only on $W$, such that $d'$ is an 
isomorphism preserving the classes of simple modules, whenever $\ell>N$. 
Thus, $\cH_k$ and $\cH_K^{(e)}$ have the same decomposition numbers, 
whenever $\ell>N$. The assertion then follows from the existence of a
canonical basic set for $\cH_K^{(e)}$, by Theorem~\ref{thmD}.
\end{proof}

\begin{thm} \label{thmC} Assume that $\ell$ is good for $G^F$ if $F$ is 
the identity on $W$, and $\ell$ does not divide the order of $W^F$, otherwise. 
Then $\cH_K^{(e)}$ and $\cH_k$ have the same number of simple modules (up 
to isomorphism). In particular, this number only depends on~$e$.
\end{thm}

Note that, even if $\cH_k$ and $\cH_K^{(e)}$ have the same number of
simple modules, $d'$ may still have a non-trivial decomposition matrix.
For examples, see the ``adjustment matrices'' computed by James
\cite{Ja}.

\begin{proof} In \cite{GeRo1}, it is shown by a general argument that the 
result is true if $\ell$ does not divide the order of $W$. Hence, for each 
type of $W$, there is only a finite number of additional cases to be 
considered.

If $W$ is of exceptional type, these additional cases can be handled by 
explicit computations with the character tables of $\cH_{K(v)}$; in the 
case where $F$ is the identity on $W$, this has been worked out explicitly 
in \cite{gennum}. Similar methods apply to the cases where $W$ is of type 
$D_4$ or $E_6$, and $W^F$ is of type $G_2$ or $F_4$, respectively. The 
desired equality can also be seen directly from the decomposition numbers 
computed in \cite[Satz~3.13.1]{myphd} (for type $G_2$ with parameters 
$q^3,q$) and Bremke \cite{Br} (for type $F_4$ with parameters $q^2,q^2,q,q$). 
In both of these cases, it is sufficient to assume that $\ell\neq 2,3$.

It remains to consider $W$ of type $A_{n-1}$, $B_n$ or $D_n$. In type
$A_{n-1}$, the result is known by Dipper--James \cite{DJ0} (the required
number is the number of $e$-regular partitions of $n$). For type $D_n$, we 
have a reduction to type $B_n$ by \cite[Theorem~6.3]{my00}. Thus, finally, 
we have to deal with the case where $F$ acts on $W$ (of type $A_{n-1}$, 
$B_n$ or $D_n$) such that $W^F$ is of type $B_m$ (where $m$ depends on $n$ 
and $F$). The weight function $L$ on $W^F$ is specified by two integers 
$a,b\geq 0$ such that:
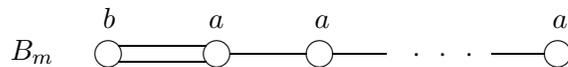
\begin{center}
\begin{picture}(250,20)
\put(  3, 03){$B_m$}
\put( 40, 05){\circle{10}}
\put( 44, 02){\line(1,0){33}}
\put( 44, 08){\line(1,0){33}}
\put( 81, 05){\circle{10}}
\put( 86, 05){\line(1,0){29}}
\put(120, 05){\circle{10}}
\put(125, 05){\line(1,0){20}}
\put(155, 02){$\cdot$}
\put(165, 02){$\cdot$}
\put(175, 02){$\cdot$}
\put(185, 05){\line(1,0){20}}
\put(210, 05){\circle{10}}
\put( 38, 15){$b$}
\put( 78, 15){$a$}
\put(118, 15){$a$}
\put(208, 15){$a$}
\end{picture}
\end{center}
We have the following possibilities:
\begin{itemize}
\item[(i)] $a=b=1$, where $G^F$ is of type $B_n$ or $C_n$ (and $m=n$).
\item[(ii)] $a=2$, $b=2s+1$,  where $G^F$ is of type ${^2\!}A_{n-1}$
(and $n=2m+s$, $s\in \{0,1\}$). 
\item[(iii)] $a=1$, $b=2$, where $G^F$ is of type ${^2\!D}_{n}$ (and 
$n=m+1$). 
\item[(iv)] $a=1$, $b=0$, where $G^F$ is of type $D_n$ (and $n=m$).
\end{itemize}
If $q^a\equiv 1\bmod \ell$, the assertion follows from
Dipper--James--Murphy \cite[Theorem 7.3]{DJM3}. (Note that, in case (ii), we
assume that $\ell$ does not divide the order of $W^F$.) Assume from now on
that $q^a\not\equiv 1 \bmod \ell$. Then, by Ariki--Mathas 
\cite[Theorem~A]{ArMa}, we have:
\begin{itemize}
\item[($*$)] {\em The number of simple $\cH_k$-modules (up to 
isomorphism) only depends on $e'$ and $A$},
\end{itemize}
where 
\begin{align*}
e'&:=\min\{j \geq 2\mid 1+q^a+q^{2a}+\cdots+q^{a(j-1)} \equiv 0 \bmod 
\ell\},\\ A&:=\{ j \in \Z \mid \theta(q^b+q^{aj})=0\}.
\end{align*}
Similarly, the number of simple $\cH_K^{(e)}$-modules only depends on $e'$
and $A_0:= \{j\in\Z\mid\zeta_e^b+\zeta_e^{aj}=0\}$. Hence, all we need
to do is to check that $A=A_0$. Now, since $q^a\not \equiv 1 \bmod \ell$,
we have $e=e'$ unless $a=2$ and $e$ is even, in which case $e'=e/2$.
Furthermore, $e$ is the multiplicative order of $q$ modulo $\ell$.  In
particular, $e$ is coprime to $\ell$. Using the fact that the canonical
map $\cO \rightarrow k$ induces an isomorphism on roots of unity of order
prime to $\ell$, it readily follows that $A=A_0$ in each of the above cases,
as required.
\end{proof}

\subsection{Proofs of Theorem~\ref{thmA} and \ref{thmB}} \label{sub-fin}
By Theorem~\ref{gero} and Corollary~\ref{lem4}, the algebra $\cH_k$
admits a canonical basic set, say $\cB_{k,q}$, if one of the following holds:
\begin{itemize}
\item $\ell$ is good for $G^F$ and $F$ acts as the identity on $W$;
\item $\ell$ is sufficiently large (where the bound comes from the proof
of Corollary~\ref{lem4}).
\end{itemize}
In these cases, Theorem~\ref{thmA} and Corollary~\ref{cons1} hold
for $G^F$ by Lemma~\ref{lem1}. 

As far as Theorem~\ref{thmB} is concerned, we first note that 
Lemma~\ref{lem1} provides a reduction to an analogous statement for 
$\cB_{k,q}$. Thus, it remains to show that $\cB_{k,q}$ only depends on $e$. 
But this follows from Theorem~\ref{thmC}, Theorem~\ref{thmD} and  
Lemma~\ref{lem3}.

%%%%%%%%%%%%%%%%%%%%%%%%%%%%%%%%%%%%%%%%%%%%%%%%%%%%%%%%%%%%%%%%%%%%%%%%%%%%%


\begin{thebibliography}{131}

\bibitem{Ar1}
{\sc S.~Ariki}, On the decomposition numbers of the Hecke algebra of
$G(m,1,n)$, J. Math. Kyoto Univ. {\bf 36} (1996), 789--808.

\bibitem{Ar2}
{\sc S.~Ariki}, On the classification of simple modules for cyclotomic
Hecke algebras of type $G(m,1,n)$ and Kleshchev multipartitions, Osaka
J. Math. {\bf 38} (2001), 827--837.

\bibitem{Ar3}
{\sc S.~Ariki}, Representations of quantum algebras and combinatorics of
Young tableaux, University Lecture Series {\bf 26},  Amer. Math. Soc.,
Providence, RI, 2002.

\bibitem{ArMa}
{\sc S.~Ariki and A.~Mathas}, The number of simple modules of the Hecke
algebras of type $G(r,1,n)$,  Math. Z. {\bf 233} (2000), 601--623.

\bibitem{Br}
{\sc K.~Bremke}, The decomposition numbers of Hecke algebras of type $F_4$
with unequal parameters, Manuscripta Math. {\bf 83} (1994), 331-346.

\bibitem{Ca2}
{\sc R.~W.~Carter}, {Finite groups of Lie type: Conjugacy classes and
complex characters}, Wiley, New York, 1985.

\bibitem{CR2}
{\sc C.~W.~Curtis and I.~Reiner}, Methods of representation theory
Vol.~I and II, Wiley, New~York, 1981 and  1987.

\bibitem{QuotHom}
{\sc R.~Dipper}, On quotients of Hom-functors and representations of
finite general linear groups I, {J.~Algebra} {\bf 130} (1990), 235--259;
II, J. Algebra {\bf 209} (1998), 199--269.

\bibitem{DJ0}
{\sc R.~Dipper and G.~D.~James}, Representations of Hecke algebras of
general linear groups, Proc.\ London Math.\ Soc. {\bf 52} (1986), 20--52.

\bibitem{DJ1}
{\sc R.~Dipper and G.D.~James}, The $q$-Schur algebra, {Proc.\ London 
Math.\ Soc.} {\bf 59} (1989), 23--50.

\bibitem{DJM3}
{\sc R.~Dipper, G. D.~James and G. E.~Murphy}, Hecke algebras of type $B_n$
at roots of unity, Proc. London Math. Soc. {\bf 70} (1995), 505--528.

\bibitem{FLOTW}
{\sc O.~Foda, B.~Leclerc, M.~Okado, J.-Y. Thibon and T.~Welsh}, Branching
functions of $A_{n-1}^{(1)}$ and Jantzen-Seitz problem for Ariki-Koike
algebras, Advances in Math. {\bf 141} (1999), 322--365.

\bibitem{fs1}
{\sc P.~Fong and B.~Srinivasan}, Brauer trees in classical groups, 
J. Algebra {\bf 131} (1990), 179--225.

\bibitem{myphd} 
{\sc M.~Geck}, Verallgemeinerte Gelfand-Graev Charaktere und 
Zer\-le\-gungs\-zah\-len end\-li\-cher Gruppen vom Lie-Typ, 
Dissertation, RWTH Aachen, 1990.

\bibitem{myuni} 
{\sc M.~Geck}, On the decomposition numbers of the finite unitary groups 
in non-defining characteristic, {Math.\ Z.} {\bf 207} (1991), 83--89.

\bibitem{mybaum}
{\sc M.~Geck}, Brauer trees of Hecke algebras, Comm.\ Algebra {\bf 20}
(1992), 2937--2973.

\bibitem{mye6}
{\sc M.~Geck}, The decomposition numbers of the Hecke algebra of type
$E_6$, Math.\ Comp.\ {\bf 61} (1993), 889--899.

\bibitem{habil}
{\sc M.~Geck}, Beitr\"age  zur Darstellungstheorie von
Iwahori--Hecke--Algebren, Habilitationsschrift, Aachener Beitr\"age zur
Mathematik {\bf 11}, Verlag der Au\-gustinus Buch\-hand\-lung, Aachen,
1995.

\bibitem{mykl}
{\sc M.~Geck}, Kazhdan-Lusztig cells and decomposition numbers. Represent.\
Theory {\bf 2} (1998), 264--277 (electronic).

\bibitem{mybourb}
{\sc M.~Geck}, Representations of Hecke algebras at roots of unity,
S\'{e}minaire Bourbaki, ann\'ee 1997-98, Ast\'erisque No.~252
(1998), Exp.~836, 33--55.

\bibitem{my00}
{\sc M.~Geck}, On the representation theory of Iwahori-Hecke algebras of
extended finite Weyl groups. Represent. Theory {\bf 4} (2000), 370--397
(electronic).

\bibitem{gennum}
{\sc M.~Geck}, On the number of simple modules of Iwahori--Hecke algebras
of finite Weyl groups. Bul. Stiit. Univ. Baia Mare, Ser. B {\bf 16}
(2000), 235--246; preprint available at 
{\tt http://arXiv.org/math.RT/0405555}

\bibitem{mylaus}
{\sc M.~Geck}, Modular representations of Hecke algebras. {\it In:}
Group representation theory, Presses Polytechniques et Universitaires
Romandes, EPFL-Press, 2006 (to appear).


\bibitem{bs1}
{\sc M.~Geck and G.~Hiss}, Basic sets of Brauer characters of finite
groups of Lie type, J.\ reine angew.\ Math. {\bf 418} (1991), 173--188.

\bibitem{lymgh}
{\sc M.~Geck and G.~Hiss}, Modular representations of finite groups of
Lie type in non-defining characteristic, {in}: Finite reductive groups:
Related structures and representations (ed.\ M.~Cabanes), pp.~195--249.
Birkh\"auser, Basel, 1997.

\bibitem{ghm1}
{\sc M.~Geck, G.~Hiss, and G.~Malle}, Cuspidal unipotent Brauer characters,
J. Algebra {\bf 168} (1994), 182--220.

\bibitem{ghm2}
{\sc M.~Geck, G.~Hiss, and G.~Malle}, Towards a classification of the
irreducible representations in non-defining characteristic of a finite
group of Lie type, Math.\ Z. {\bf 221} (1996), 353--386.

\bibitem{GeJa}
{\sc M.~Geck and N.~Jacon}, Canonical basic sets in type $B$, J. Algebra
({\em to appear}).

\bibitem{GeLu}
{\sc M.~Geck and K.~Lux}, The decomposition numbers of the Hecke algebra of
type $F_4$. Manus\-cripta Math.\ {\bf 70} (1991), 285--306.

\bibitem{GeMa}
{\sc M.~Geck and G.~Malle}, On the existence of a unipotent support for the
irreducible characters  of finite groups of Lie type.  Trans. Amer. Math.
Soc. {\bf 352} (2000), 429--456.

\bibitem{ourbuch}
{\sc M.~Geck and G.~Pfeiffer}, Characters of finite Coxeter groups and
Iwahori--Hecke algebras, London Math. Soc. Monographs, New Series {\bf 21},
Oxford University Press, New York 2000. xvi+446 pp.

\bibitem{GeRo1}
{\sc M.~Geck and R.~Rouquier}, Centers and simple modules for Iwahori-Hecke
algebras.  {\it In}: Finite reductive groups: Related structures and
representations (ed.\ M.~Cabanes), pp.~251--272.  Birkh\"auser, Basel, 1997.

\bibitem{GeRo2}
{\sc M.~Geck and R.~Rouquier}, Filtrations on projective modules
for Iwahori--Hecke algebras. {\it In}: Modular Representation
Theory of Finite Groups (Charlottesville, VA, 1998; eds. M.~J.~Collins,
B.~J.~Parshall and L.~L.~Scott), p.~211--221, Walter de Gruyter, Berlin 2001.

\bibitem{H2} 
{\sc G.~Hiss}, Decomposition numbers of finite groups of Lie type in 
non-defining characteristic, {\em in:} G.~O.~Michler and C.~M.~Ringel, 
Eds., Representation Theory of Finite Groups and Finite-Dimensional 
Algebras, Birkh{\"a}user, 1991, pp.~405--418.

\bibitem{Hiss1}
{\sc G.~Hiss}, Harish-Chandra series of Brauer characters in a finite group 
with a split $BN$-pair, J.~London Math.\ Soc. {\bf 48} (1993), 219--228.

\bibitem{Jac0}
{\sc N.~Jacon}, Repr\'esentations modulaires des alg\`ebres de Hecke et
des alg\`ebres de Ariki-Koike, Th\`ese de Doctorat, Universit\'e 
Lyon 1, 2004; available at ``theses-ON-line''
{\tt http://tel.ccsd.cnrs.fr/documents/archives0/00/00/63/83}.

\bibitem{Jac1}
{\sc N.~Jacon}, Sur les repr\'esentations modulaires des alg\`ebres de
Hecke de type $D\sb n$, J. Algebra  {\bf 274}  (2004), 607--628

\bibitem{Jac2}
{\sc N.~Jacon}, On the parametrization of the simple modules for
Ariki-Koike algebras at roots of unity,  J. Math. Kyoto Univ. {\bf 44}
(2004),  729--767

\bibitem{Ja}
{\sc G.~D.~James}, The decomposition matrices of $\GL_n(q)$ for
$n \leq 10$, Proc.\ London Math.\ Soc. {\bf 60} (1990), 225--265.

\bibitem{Kaw1}
{\sc N.~Kawanaka}, Generalized Gelfand-Graev representations and
Ennola Duality. {\em In:} Algebraic Groups and Related Topics, Advanced 
Studies in Pure Math.\ {\bf 6}, Kinokuniya, Tokyo, and North-Holland, 
Amsterdam, 1985, pp.~175--206.

\bibitem{L1} 
{\sc G.~Lusztig}, Irreducible representations of finite classical groups, 
Invent.\ Math. {\bf 43} (1977), 125--175.

\bibitem{LuBC}
{\sc G.~Lusztig}, {On a theorem of {B}enson and {C}urtis},
J. Algebra {\bf 71} (1981), 490--498.

\bibitem{Lu4}
{\sc G.~Lusztig}, Characters of reductive groups over a finite field,
Annals Math.\ Studies, vol. 107, Princeton University Press, 1984.

\bibitem{L6} 
{\sc G.~Lusztig}, A unipotent support for irreducible
representations, Adv.\ Math.\ {\bf 94} (1992), 139--179.

\bibitem{Lusztig03}
{\sc G.~Lusztig}, Hecke algebras with unequal parameters, CRM Monographs
Ser.~{\bf 18}, Amer. Math. Soc., Providence, RI, 2003.

\bibitem{muell}
{\sc J.~M\"uller}, {Zerlegungszahlen f\"ur generische Iwahori-Hecke-Algebren
von exzeptionellem Typ}, Dissertation, RWTH Aachen, 1995.

\end{thebibliography}
\end{document}